\documentclass[12pt,a4paper]{amsart}

\usepackage{amscd,amssymb,amsthm}
\usepackage{graphicx}

\theoremstyle{plain}
\newtheorem{theorem}{Theorem}[section]
\newtheorem{corollary}[theorem]{Corollary}
\newtheorem{lemma}[theorem]{Lemma}

\theoremstyle{definition}

\theoremstyle{remark}
\newtheorem{remark}{Remark}[section]

\numberwithin{equation}{section}

\numberwithin{table}{section}

\numberwithin{figure}{section}

\def\dj{d{\hspace{-.32em}\raisebox{.7ex}{-}}}
\def\Dj{D{\hspace{-.75em}\raisebox{.3ex}{-}\hspace{.4em}}}

\title[Some inequalities for Kurepa's function]{SOME INEQUALITIES FOR KUREPA'S FUNCTION}
\author{Branko J. Male\v sevi\' c}
\address{University of Belgrade    \\
Faculty of Electrical Engineering  \\
Serbia \& Montenegro}              \email{malesevic@etf.bg.ac.yu}

\keywords{Kurepa's function, Inequalities for integrals} \subjclass[2000]{26D15}

\begin{document}

\begin{abstract}
In this paper we consider Kurepa's function $K(z)$
\cite{Kurepa_71}. We give some recurrent relations for Kurepa's
function via appropriate sequences of rational functions and gamma
function. Also, we give some inequalities for Kurepa's function
$K(x)$ for positive values of $x$.
\end{abstract}

\maketitle

\footnotetext{Research partially suported by the MNTRS,
Serbia \& Montenegro, Grant No. 1861.}

\section{Kurepa's function $K(z)$}\label{sec_1}

\Dj uro Kurepa considered, in article \cite{Kurepa_71}, the
function of left factorial $!n$ as a sum of factorials $!n = 0! +
1! + 2! + \ldots + (n\!-\!1)!$. Let us use the standard notation:
\begin{equation}
\label{K_SUM_1}
K(n) = \displaystyle\sum\limits_{i=0}^{n-1}{i!}.
\end{equation}
Sum (\ref{K_SUM_1}) corresponds to the sequence $A003422$ in
\cite{Sloane_03}. An analytical extension of the function
(\ref{K_SUM_1}) over the set of complex numbers is determined by
the integral:
\begin{equation}
\label{K_INT_1}
K(z)
=
\displaystyle\int\limits_{0}^{\infty}{
e^{-t} \displaystyle\frac{t^{z}-1}{t-1} \: dt},
\end{equation}
which converges for $\mbox{Re} \: z > 0$~\cite{Kurepa_73}. For function $K(z)$ we use
the term {\em Kurepa's function}. It is easily verified that Kurepa's function $K(z)$
is a solution of the functional equation:
\begin{equation}
\label{K_FE_1}
K(z) - K(z-1) = \Gamma(z).
\end{equation}
Let us observe that since $K(z-1) = K(z) - \Gamma(z)$, it is possible to make
the analytic continuation of Kurepa's function $K(z)$ for \mbox{$\mbox{Re} \, z \leq 0$}.
In that way, the Kurepa's function $K(z)$ is a meromorphic function with simple
poles at $z = -1$~and~$z =-n$ \mbox{$(n\!\geq\!3)$}~\cite{Kurepa_73}.
Let us emphasize that in the following consideration, in sections \ref{sec_2}
and \ref{sec_3}, it is sufficient to use only the fact that function $K(z)$ is
a solution of the functional equation~(\ref{K_FE_1}).

\section{Representation of the Kurepa's function via sequences of polynomials
and gamma function}\label{sec_2}

\Dj uro Kurepa considered, in article \cite{Kurepa_73}, the sequences
of following polynomials:
\begin{equation}
P_{n}(z) = (z-n)P_{n-1}(z) + 1,
\end{equation}
with an initial member $P_{0}(z)=1$. On the basis of
\cite{Kurepa_73} we can conclude that the following statements are
true:
\begin{lemma}
For each $n \!\in\! N$ and $z \!\in\! \mbox{\bf C}$ we have explicitly$:$
\begin{equation}
P_{n}(z)
=
1+\displaystyle\sum\limits_{j=0}^{n-1}{\displaystyle\prod\limits_{i=0}^{j}{(z-n+i)}}.
\end{equation}
\end{lemma}
\begin{theorem}
For each $n \!\in\! N$ and $z \!\in\! \mbox{\bf C} \backslash
\mbox{\big (}\mathop{\mbox{\bf Z}}^{-} \cup \, \{0,1,\ldots,n\} \mbox{\big )}$ is valid$:$
\begin{equation}
\label{K_P_Veza}
K(z) = K(z - n) + {\big (} P_{n}(z) - 1 {\big )} \cdot \Gamma(z-n).
\end{equation}
\end{theorem}

\section{Representation of the Kurepa's function via sequences of
rational functions and gamma function}\label{sec_3}

Let us observe that on the basis of a functional equation for
gamma function \mbox{$\Gamma(z+1)$} \mbox{$ = z\Gamma(z)$}, it
follows that the Kurepa's function is the solution of the
following functional equation:
\begin{equation}
\label{K_FE_2}
K(z+1) - (z+1)K(z) + zK(z-1) = 0.
\end{equation}
For $z \!\in\! \mbox{\bf C} \backslash \{0\}$, based on (\ref{K_FE_2}), it  is valid:
\begin{equation}
\label{K_FE_2_1}
K(z-1)
=
\mbox{\small $\displaystyle\frac{z+1}{z}K(z)$}
-
\mbox{\small $\displaystyle\frac{1}{z}K(z+1)$}
=
Q_{1}(z) K(z) - R_{1}(z) K(z+1),
\end{equation}
for rational functions $Q_{1}(z) \!=\! \mbox{\small $\displaystyle\frac{z\!+\!1}{z}$}$,
$R_{1}(z) \!=\! \mbox{\small $\displaystyle\frac{1}{z}$}$ over $\mbox{\bf C}
\backslash \{0\}$. Next,  for $z \!\in\! \mbox{\bf C} \backslash \{0,1\}$,
based on (\ref{K_FE_2}), it is valid:
\begin{equation}
\begin{array}{rcl}
K(z-2)
&\!\!\!\!=\!\!\!\!&
\mbox{\small $\displaystyle\frac{z}{z-1}$}K(z-1)
-
\mbox{\small $\displaystyle\frac{1}{z-1}$}K(z)                                \\[0.5 ex]
&\!\!\!\!\mathop{=}\limits_{(\ref{K_FE_2_1})}\!\!\!\!&
\mbox{\small $\displaystyle\frac{z}{z-1}$}{\Big (}
\mbox{\small $\displaystyle\frac{z+1}{z}$}K(z)
-
\mbox{\small $\displaystyle\frac{1}{z}$}K(z+1){\Big )}
-
\mbox{\small $\displaystyle\frac{1}{z-1}$}K(z)                                \\[0.5 ex]
&\!\!\!\!=\!\!\!\!&
\mbox{\small $\displaystyle\frac{z}{z-1}K(z)$}
-
\mbox{\small $\displaystyle\frac{1}{z-1}K(z+1)$}
=
Q_{2}(z) K(z) - R_{2}(z) K(z+1),
\end{array}
\end{equation}
for rational functions $Q_{2}(z) \!=\! \mbox{\small $\displaystyle\frac{z}{z\!-\!1}$}$,
$R_{2}(z) \!=\! \mbox{\small $\displaystyle\frac{1}{z\!-\!1}$}$ over $\mbox{\bf C} \backslash
\{0,1\}$. Thus, for values $z \!\in\! \mbox{\bf C} \backslash \{0,1,\ldots,n\!-\!1\}$,
based on (\ref{K_FE_2}), by mathematical induction it is true:
\begin{equation}
K(z-n)
=
Q_{n}(z) K(z) - R_{n}(z) K(z+1),
\end{equation}
for rational functions $Q_{n}(z)$, $R_{n}(z)$ over $\mbox{\bf C}
\backslash \{0,1,\ldots,n\!-\!1\}$,  which fulfill the same
recurrent relations:
\begin{equation}
\label{Rec_Q}
Q_{n}(z)
=
\mbox{\small $\displaystyle\frac{z-n+2}{z-n+1}$}Q_{n-1}(z)
-
\mbox{\small $\displaystyle\frac{1}{z-n+1}$}Q_{n-2}(z)
\end{equation}
and
\begin{equation}
\label{Rec_R}
R_{n}(z)
=
\mbox{\small $\displaystyle\frac{z-n+2}{z-n+1}$}R_{n-1}(z)
-
\mbox{\small $\displaystyle\frac{1}{z-n+1}$}R_{n-2}(z),
\end{equation}
with  different initial functions $Q_{1,2}(z)$ and $R_{1,2}(z)$.

\smallskip
Based on the previous consideration we can conclude:
\begin{lemma}
For each $n \!\in\! N$ and $z \!\in\! \mbox{\bf C} \backslash \{0,1,\ldots,n\!-\!1\}$
let the rational function $Q_{n}(z)$ be determined by the recurrent relation
{\rm (\ref{Rec_Q})} with initial functions
$Q_{1}(z) = \mbox{\small $\displaystyle\frac{z+1}{z}$}$
and
$Q_{2}(z) = \mbox{\small $\displaystyle\frac{z}{z-1}$}$.
Thus the sequence $Q_{n}(z)$ has an explicit form$:$
\begin{equation}
Q_{n}(z)
=
1
+
\displaystyle\sum\limits_{j=0}^{n-1}{\:
\displaystyle\prod\limits_{i=0}^{j}{
\displaystyle\frac{1}{z-i}}}.
\end{equation}
\end{lemma}
\begin{lemma}
For each $n \!\in\! N$ and $z \!\in\! \mbox{\bf C} \backslash \{0,1,\ldots,n\!-\!1\}$
let the rational function $R_{n}(z)$ be determined by the recurrent relation
{\rm (\ref{Rec_R})} with initial functions
$R_{1}(z) = \mbox{\small $\displaystyle\frac{1}{z}$}$
and
$R_{2}(z) = \mbox{\small $\displaystyle\frac{1}{z-1}$}$.
Thus the sequence $R_{n}(z)$ has an explicit form$:$
\begin{equation}
R_{n}(z)
=
\displaystyle\sum\limits_{j=0}^{n-1}{\:
\displaystyle\prod\limits_{i=0}^{j}{
\displaystyle\frac{1}{z-i}}}.
\end{equation}
\end{lemma}
\begin{theorem}
For each $n \!\in\! N$ and $z \!\in\! \mbox{\bf C} \backslash
\mbox{\big (}\mathop{\mbox{\bf Z}}^{-} \cup \, \{0,1,\ldots,n\!-\!1\}
\mbox{\big )}$ it is valid$:$
\begin{equation}
\label{K_Q_Veza}
K(z) = K(z\!-\!n) + {\big (} Q_{n}(z) - 1 {\big )} \cdot \Gamma(z+1)
\end{equation}
and
\begin{equation}
\label{K_R_Veza}
K(z) = K(z\!-\!n) + R_{n}(z) \cdot \Gamma(z+1).
\end{equation}
\end{theorem}

\section{Some inequalities for Kurepa's function}\label{sec_4}

In this section we consider the Kurepa's function $K(x)$, given by an integral representation
(\ref{K_INT_1}), for positive values of $x$. Thus the Kurepa's function is positive and
in the following consideration we give some inequalities for the Kurepa's function.
\begin{lemma}
\label{lm_11}
For $x \! \in \! [0,1]$ the following inequalities are true$:$
\begin{equation}
\label{lm_11_iq_1}
\Gamma(x+1/2) < x^2 - \displaystyle\frac{7}{4}x + \displaystyle\frac{9}{5}
\end{equation}
and
\begin{equation}
\label{lm_11_iq_2}
(x+2) \Gamma(x+1) > \displaystyle\frac{9}{5}.
\end{equation}
\end{lemma}
\begin{proof}
It is sufficient to use an approximation formula for the function
$\Gamma(x+1)$ with polynomial of the fifth degree:
$P_{5}(x)=-0.1010678\,x^5+0.4245549\,x^4-0.6998588\,x^3
+0.9512363\,x^2-0.5748646\,x+1$ which has an absolute error
$|\varepsilon(x)| < 5 \cdot 10^{-5}$ for values of argument $x
\!\in\! [0,1]$ \cite{Abramowitz_&_Stegun_72} (formula 6.1.35., page 257.).
To prove the first inequality, for values $x \!\in\! [0,1/2]$, it is
necessary to consider an equivalent  inequality obtained by the
following substitution $t = x + 1/2$ {\big (}thus $\Gamma(x+1/2) =
\Gamma(t+1)/t${\big )}. To prove the first inequality, for values
$x \!\in\! (1/2,1]$, it is necessary to consider an equivalent
inequality by the following substitution $t = x - 1/2$ {\big
(}thus $\Gamma(x+1/2) = \Gamma(t+1)${\big )}.
\end{proof}
\begin{remark}
Let us notice that for a proof  of the previous inequalities it is
possible to use other polynomial approximations {\rm (}of a lower
degree{\rm )} of functions $\Gamma(x+1/2)$ and $\Gamma(x+1)$ for
values $x \!\in\! [0,1]$.
\end{remark}
\begin{lemma}
\label{lm_12}
For $x \! \in \! [0,1]$ the following inequality is true$:$
\begin{equation}
\label{lm_12_iq_1}
K(x) \leq \displaystyle\frac{9}{5}x.
\end{equation}
\end{lemma}
\begin{proof}
Let us notice that the first derivation of the Kurepa's function $K(x)$,
for values $x \! \in \! [0,1]$, is given by the following integral \cite{Kurepa_73}:
\begin{equation}
K^{'}(x)= \displaystyle\int\limits_{0}^{\infty}{
e^{-t} t^x \displaystyle\frac{\log t}{t-1} \: dt}.
\end{equation}
For $t \! \in \! (0, \infty) \backslash \{1\}$ Karamata's
inequality is true: $\displaystyle\frac{\log t}{t-1} \leq
\displaystyle\frac{1}{\sqrt{t}}$ \cite{Karamata_60}. Hence, for $x
\! \in \! [0,1]$ the following inequality is true:
\begin{equation}
\label{lm_12_iq_2}
K^{'}(x)
=
\displaystyle\int\limits_{0}^{\infty}{e^{-t} t^x
\displaystyle\frac{\log t}{t-1} \: dt}
\leq
\displaystyle\int\limits_{0}^{\infty}{e^{-t} t^{x-1/2} \: dt}
=
\Gamma(x+1/2).
\end{equation}
Next, on the basis of the Lemma \ref{lm_11} and inequality (\ref{lm_12_iq_2}),
for $x \! \in \! [0,1]$, the following inequalities are true:
\begin{equation}
K(x)
\leq
\displaystyle\int\limits_{0}^{x}{\Gamma(t+1/2)\:dt}
\leq
\displaystyle\int\limits_{0}^{x}{{\Big (}
t^2\!-\!\displaystyle\frac{7}{4}t\!+\!\displaystyle\frac{9}{5}{\Big )} \: dt}
\leq
\displaystyle\frac{9}{5}x.
\end{equation}
\end{proof}
\begin{theorem}
\label{lm_13}
For $x \geq 3$ the following inequality is true$:$
\begin{equation}
\label{Ineq_01}
K(x-1) \leq \Gamma(x),
\end{equation}
while the equality is true for $x = 3$.
\end{theorem}
\begin{proof}
Based on the functional equation (\ref{K_FE_1}) the inequality (\ref{Ineq_01}),
for $x \geq 3$, is equivalent to the following inequality:
\begin{equation}
\label{lm_13_iq_1}
K(x) \leq  2\Gamma(x).
\end{equation}
Let us represent $[3,\infty) =
\bigcup_{n=3}^{\infty}{[n,n\!+\!1)}$. Then, we prove that the
inequality (\ref{lm_13_iq_1}) is true, by mathematical induction
over intervals $[n,n\!+\!1)$ ($n\!\geq\!3$).

\smallskip\noindent
{\boldmath $(i)$} Let's $x \!\in\! [3,4)$. Then following decomposition is true:
$K(x) = K(x-3) + \Gamma(x-2) + \Gamma(x-1) + \Gamma(x)$. Hence, by Lemma \ref{lm_12},
the following inequality is true:
\begin{equation}
\label{lm_13_iq_2}
K(x)
\leq
\displaystyle\frac{9}{5}(x\!-\!3)\!+\!\Gamma(x\!-\!2)\!+\!\Gamma(x\!-\!1)\!+\!\Gamma(x),
\end{equation}
because $x-3 \!\in\! [0,1)$. Next, by Lemma \ref{lm_11}, the following inequality is true:
\begin{equation}
\label{lm_13_iq_3}
\displaystyle\frac{9}{5} (x\!-\!3)
\leq
(x\!-\!1)(x\!-\!3) \Gamma(x\!-\!2),
\end{equation}
because $x-3 \!\in\! [0,1)$. Now, based on (\ref{lm_13_iq_2}) and
(\ref{lm_13_iq_3}) we conclude that the inequality is true:
\begin{equation}
K(x)
\leq
(x\!-\!1)(x\!-\!3) \Gamma(x\!-\!2)
+
\Gamma(x\!-\!2)
+
\Gamma(x\!-\!1)
+
\Gamma(x)
=
2\Gamma(x).
\end{equation}

\smallskip
\noindent
{\boldmath $(ii)$} Let the inequality (\ref{lm_13_iq_1}) be true
for $x \!\in\! [n,n\!+\!1)$ $(n \!\geq\! 3)$.

\smallskip \noindent
{\boldmath $(iii)$} For $x \!\in\! [n\!+\!1,n\!+\!2)$ $(n \!\geq\! 3)$,
based on the inductive hypothesis, the following inequality is true:
\begin{equation}
K(x) = K(x-1) + \Gamma(x) \leq 2\Gamma(x-1) + \Gamma(x) \leq 2\Gamma(x).
\end{equation}
\end{proof}
\begin{remark}
The inequality {\rm (\ref{lm_13_iq_1})} is an improvement of
inequalities of Aran\dj elovi\' c$:\;$ $K(x) \leq 1 + 2\Gamma(x)$,
given in {\rm \cite{Kurepa_73}}, with respect to the
interval $[3,\infty)$.
\end{remark}
\begin{corollary}
For each $k \!\in\! N$ and $x \geq k\!+\!2$ the following inequality is true$:$
\begin{equation}
\label{Ineq_GK1}
\displaystyle\frac{K(x\!-\!k)}{\Gamma(x\!-\!k\!+\!1)} \leq 1,
\end{equation}
while the equality is true for $x = k\!+\!2$.
\end{corollary}
\begin{theorem}
For each $k \!\in\! N$ and $x \geq k\!+\!2$ the following double inequality is true$:$
\begin{equation}
\label{Ineq_GK2}
R_{k}(x)
<
\displaystyle\frac{K(x)}{\Gamma(x+1)}
\leq
\displaystyle\frac{P_{k-1}(x)+1}{P_{k-1}(x)} \cdot R_{k}(x),
\end{equation}
while the equality is true  for $x = k\!+\!2$.
\end{theorem}
\begin{proof}
For each $k \!\in\! N$ and $x\!>\!k$ let's introduce the following function
$G_{k}(x) = \sum\limits_{i=0}^{k-1}{\Gamma(x\!-\!i)}$.
Thus, the following relations:
\begin{equation}
\label{G_R_Veza}
G_{k}(x) = \Gamma(x+1) \cdot R_{k}(x)
\end{equation}
and
\begin{equation}
\label{G_P_Veza}
G_{k}(x) = \Gamma(x-k) \cdot (P_{k}(x)-1)
\end{equation}
are true. The inequality $G_{k}(x) < K(x)$ is true for $x > k$. Hence, based on
(\ref{G_R_Veza}), the left inequality in (\ref{Ineq_GK2}) is true for all $x \geq k\!+\!2$.
On the other hand, based on (\ref{G_P_Veza}) and (\ref{Ineq_GK1}), for $x \geq k\!+\!2$,
the following inequality is true:
\begin{equation}
\begin{array}{rcl}
\displaystyle\frac{K(x)}{G_{k}(x)}
& \!\!=\!\! &
1\!+\!\displaystyle\frac{K(x\!-\!k)}{G_{k}(x)}
\:=\:
1\!+\!\displaystyle\frac{K(x\!-\!k)}{\Gamma(x\!-\!k)(P_{k}(x)\!-\!1)}        \\[2.5 ex]
& \!\!=\!\! &
1\!+\!\displaystyle\frac{K(x\!-\!k)/\Gamma(x\!-\!k\!+\!1)}{P_{k-1}(x)}
\:\leq\:
1\!+\!\displaystyle\frac{1}{P_{k-1}(x)}
\:=\:
\displaystyle\frac{P_{k-1}(x)+1}{P_{k-1}(x)}.
\end{array}
\end{equation}
Hence, based on (\ref{G_R_Veza}), the right inequality in (\ref{Ineq_GK2}) is true
for all $x \geq k\!+\!2$.
\end{proof}
\begin{corollary}
If for each $k \!\in\! N$ we mark$:$
\begin{equation}
A_{k}(x) = R_{k}(x)
\quad \mbox{and} \quad
B_{k}(x) = \displaystyle\frac{P_{k-1}(x)+1}{P_{k-1}(x)} \cdot R_{k}(x),
\end{equation}
thus, the following is true$:$
\begin{equation}
A_{k}(x)
<
A_{k+1}(x)
<
\displaystyle\frac{K(x)}{\Gamma(x+1)}
\leq
B_{k+1}(x)
<
B_{k}(x)
\quad (x \geq k\!+\!3)
\end{equation}
and
\begin{equation}
A_{k}(x), B_{k}(x) \sim \displaystyle\frac{1}{x}
\;\;\wedge\;\;
B_{k}(x) - A_{k}(x) = \displaystyle\frac{R_{k}(x)}{P_{k-1}(x)}
\sim
\displaystyle\frac{1}{x^k} \quad (x \rightarrow \infty).
\end{equation}
\end{corollary}


\bigskip

\begin{thebibliography}{9}
\setlength{\itemsep}{5pt}

\bibitem{Abramowitz_&_Stegun_72}
M. ABRAMOWITZ \textsc{and} I. A. STEGUN:
Handbook of Mathematical Functions With Formulas, Graphs, and Mathematical Tables,
{\em USA National Bureau of Standards, Applied~Math.~Series~-~55, Tenth Printing}
(1972). 

\bibitem{Karamata_60}
J. KARAMATA: Sur quelques problemes poses par Ramanujan,
{\em J. Indian Math Soc.} (N.S.) {\bf 24} (1960), 343-365.

\bibitem{Kurepa_71}
\Dj . KUREPA: On the left factorial function $!n$,
{\em Mathematica Balkanica} {\bf 1} (1971), 147-153.

\bibitem{Kurepa_73}
\Dj . KUREPA:
Left factorial in complex domain,
{\em Mathematica Balkanica} {\bf 3} (1973), 297-307.

\bibitem{Sloane_03}
N. J. A. SLOANE: The-On-Line Encyclopedia of Integer Sequences.

Available online at {\tt http://www.research.att.com/\verb|~|njas/sequences/}

\end{thebibliography}
\end{document}